\documentstyle[amsfonts, txfonts, amssymb, latexsym, 12pt]{amsart}

\newtheorem{theorem}{Theorem}[section]

\newtheorem{problem}[theorem]{Problem}
\def\Z{{\mbox{\rm\kern.25em
\vrule width.03em height0.57ex depth0ex
\kern.033em
\vrule width.03em height1.52ex depth-0.96ex \kern-.338em Z}}}
\def\R{{\mbox{\rm I\kern-.22em R}}}
\def\N{{\mbox{\rm I\kern-.22em N}}}

\def\D{{\bf D}}

\def\size{{\rm size}}

\def\energy{{\rm energy}}

\def\M{{\cal{M}}}

\def\D{{\cal{D}}}
\def\J{{\cal{J}}}
\def\I{{\cal{I}}}
\def\C{{\cal{C}}}
\def\M{{\cal{M}}}
\def\F{{\cal{F}}}

\def\111{\gamma}

\def\be#1{\begin{equation}\label{#1}}
\def\bas{\begin{align*}}
\def\eas{\end{align*}}
\def\bi{\begin{itemize}}
\def\ei{\end{itemize}}
\newenvironment{proof}{\noindent {\bf Proof} }{\endprf\par}
\def \endprf{\hfill  {\vrule height6pt width6pt depth0pt}\medskip}
\def\emph#1{{\it #1}}

\title{Flag Paraproducts}

\author{Camil Muscalu}
\address{Department of Mathematics, Cornell University, Ithaca, NY 14853}
\email{camil@@math.cornell.edu}

\begin{document}

\begin{abstract}
We describe the theory of flag paraproducts and their relationship to the field of differential equations.
\end{abstract}

\maketitle

\section{Short Introduction}

The main goal of the present paper is to describe the theory of a new class of 
multi-linear operators which we named ``paraproducts with flag singularities'' 
(or in short ``flag paraproducts''). 

These objects, which have been introduced in \cite{c} as being generalizations of the
``lacunary versions`` of the ``bi-est operators`` of \cite{mtt:walshbiest},
\cite{mtt:fourierbiest}, \cite{mtt:multiest}, turned out in the meantime to have very natural
connections to several problems in the theory of differential equations.

While most of the article is expository, we also prove as a consequence of our discussion
a new ``paradifferential Leibnitz rule'', which may be of independent interest.

In Section 2 we briefly recall the theory of classical paraproducts and then, in Section 3,
we present the basic facts about the flag paraproducts. Sections 4, 5 and 6 are devoted to the
description of the various connections of the flag paraproducts: first, to the AKNS systems
of mathematical physics and scattering theory, then to what we called ``the grand Leibnitz rule''
for generic non-linearities and in the end to the theory of non-linear Schr\"{o}dinger
equations. The last section, Section 7, presents a sketch of some of the main ideas needed to
understand the boundedness properties of these flag paraproducts.

{\bf Acknowledgements}:
The present article is based on the author's lecture at the ``8th International Conference
on Harmonic Analysis and PDE'' held in El Escorial - Madrid, in June 2008.
We take this opportunity to thank the organizers once more for the invitation and for their 
warm hospitality during our stay in Spain.
We are also grateful to the NSF for partially supporting this work.

\section{Classical Paraproducts}

If $n\geq 1$, let us denote by $T$ the $n$ - linear singular integral operator given by

\begin{equation}\label{1}
T(f_1, ..., f_n)(x) = \int_{\R^n}
f_1(x-t_1) ... f_n(x-t_n) K(t) d t,
\end{equation}
where $K$ is a Calder\'{o}n - Zygmund kernel \cite{stein}.

Alternatively, $T$ can also be written as

\begin{equation}\label{2}
T_m(f_1, ..., f_n)(x) = \int_{\R^n} m(\xi)
\widehat{f_1}(\xi_1) ... 
\widehat{f_n}(\xi_n) 
e^{ 2\pi i x (\xi_1 + ... + \xi_n)}
d \xi,
\end{equation}
where $m(\xi) = \widehat{K}(\xi)$ is a classical multiplier, satisfying the well known 
Marcinkiewicz - Mihlin - H\"{o}rmander condition

\begin{equation}\label{3}
|\partial^{\alpha} m(\xi)| \lesssim \frac{1}{|\xi|^{|\alpha|}}
\end{equation}
for sufficiently many multi-indices $\alpha$.
\footnote{We use the standard notation $A\lesssim B$ to denote the fact that there exists
a constant $C > 0$ so that $A \leq C\cdot B$. We also denote by $\M(\R^n)$ the class of all
such multipliers.}

These operators play a fundamental role in analysis and PDEs and they are called 
``paraproducts''. \footnote{It is easy to observe that in the particular case when
$m=1$, $T_m(f_1, ..., f_n)(x)$ becomes the product of the $n$ functions 
$f_1(x)\cdot ... \cdot f_n(x)$. Also, as stated, the formulas are for functions of one variable, but the whole theory extends easily to an arbitrary euclidean
space $\R^d$.}
The following Coifman - Meyer theorem is a classical result in harmonic analysis
\cite{meyerc}, \cite{ks}, \cite{gt}.

\begin{theorem}
For every $m\in \M(\R^n)$, the $n$-linear multiplier $T_m$ maps
$L^{p_1}\times ... \times L^{p_n} \rightarrow L^p$ boundedly, as long as
$1< p_1, ..., p_n \leq \infty$, $1/p_1 + ... + 1/p_n = 1/p$ and $0<p<\infty$.
\end{theorem}

To recall some of the main ideas which appear in the proof of the theorem, let us
assume that the kernel $K(t)$ has the form

\begin{equation}\label{4}
K(t) = \sum_{k\in \Z} \Phi_k^1(t_1) ... \Phi_k^n(t_n),
\end{equation}
where each $\Phi_k^j$ is an $L^1$ - normalized bump function adapted to the interval
$[-2^{-k}, 2^{-k}]$
\footnote{In fact, modulo some technical issues, one can always assume that this is the case.}.

As a consequence, for any $1<p<\infty$, one has

$$\left\|T_m(f_1, ..., f_n)\right\|_p =
\left| \int_{\R}T_m(f_1, ..., f_n)(x) f_{n+1}(x) d x \right| =
$$

\begin{equation}\label{5}
\left|
\int_{\R}\sum_{k\in\Z}
(f_1\ast\Phi_k^1)(x) ...
(f_n\ast\Phi_k^n)(x)
(f_{n+1}\ast\Phi_k^{n+1})(x)
d x
\right|,
\end{equation}
where $f_{n+1}$ is a well chosen function with $\|f_{n+1}\|_{p'} = 1$ (and $1/p+1/p'=1$), while
the family $(\Phi_k^{n+1})_k$ is also as usual well chosen so that the above equality 
holds true. One should also recall the standard fact that since $K$ is a 
Calder\'{o}n - Zygmund kernel, one can always assume that at least two of the families
$(\Phi_k^j)_k$ for $j=1, ..., n+1$ are of ``$\Psi$ type'', in the sense that the Fourier
transform of the corresponding $k$th term is supported in 
$[-2^{k+1}, -2^{k-1}]\cup [2^{k-1}, 2^{k+1}]$, while all the others are of ``$\Phi$ type'',
in the sense that the Fourier transform of the corresponding $k$th term is supported
in $[-2^{k+1}, 2^{k+1}]$. For simplicity, we assume that in our case
$(\Phi_k^1)_k$ and $(\Phi_k^2)_k$ are of ``$\Psi$ type''
\footnote{We will use this ``$\Psi$ - $\Phi$'' terminology throughout the paper.}.
 
Then, (\ref{5}) can be majorized by

$$\int_{\R}
\left(\sum_k |f_1\ast\Phi_k^1(x)|^2\right)^{1/2}
\left(\sum_k |f_2\ast\Phi_k^2(x)|^2\right)^{1/2}
\prod_{j\neq 1,2}
\sup_k |f_j\ast\Phi_k^j(x)| d x \lesssim
$$

$$\int_{\R} S f_1(x) \cdot S f_2(x) \cdot \prod_{j\neq 1,2} M f_j(x) dx$$
where $S$ is the Littlewood - Paley square function and $M$ is the Hardy - Littlewood maximal
function.

Using now their boundedness properties \cite{stein}, one can easily conclude that
$T_m$ is always bounded from $L^{p_1}\times ... \times L^{p_n} \rightarrow L^p$, as long as
all the indices $p_1, ..., p_n, p$ are strictly between 1 and $\infty$. The $L^{\infty}$
case is significantly harder and it usually follows from the so called $T1$ - theorem of
David and Journ\'{e} \cite{stein}. Once the ``Banach case'' of the theorem has been understood,
the ``quasi - Banach case'' follows from it by using Calder\'{o}n - Zygmund decompositions
for all the functions $f_1, ..., f_n$ carefully \cite{meyerc}, \cite{ks}, \cite{gt}.

\section{Flag Paraproducts}

We start with the following concrete example

\begin{equation}\label{6}
T(f,g,h)(x) = \int_{\R^7}
f(x-\alpha_1-\beta_1)
g(x-\alpha_2-\beta_2-\gamma_1)
h(x-\alpha_3-\gamma_2) K(\alpha) K(\beta) K(\gamma)
d \alpha
d \beta
d \gamma
\end{equation}
which is a prototype of a ``flag paraproduct''. As one can see, there are now
three kernels acting on our set of three functions. $K(\beta)$ and $K(\gamma)$
being kernels of two variables, act on the pairs $(f,g)$ and $(g,h)$ respectively, while
$K(\alpha)$ being a kernel of three variables acts on all three functions $(f,g,h)$ and all of 
them in a ``paraproduct manner''. The point is that all these three ``actions'' happen
simultaneously.

Alternatively, one can rewrite (\ref{6}) as

$$
T(f,g,h)(x) = \int_{\R^3} m(\xi)
\widehat{f}(\xi_1) 
\widehat{g}(\xi_2)
\widehat{h}(\xi_3)
 e^{ 2\pi i x (\xi_1 + \xi_2 +\xi_3)}
d \xi$$
where 

$$m(\xi) = m'(\xi_1, \xi_2)\cdot m''(\xi_2, \xi_3)\cdot m'''(\xi_1, \xi_2, \xi_3)$$
is now a product of three classical symbols, two of them in $\M(\R^2)$ and the third
in $\M(\R^3)$.

Generally, for $n\geq 1$, we denote by $\M_{flag}(\R^n)$ the set of all symbols $m$ given by
arbitrary products of the form

$$m(\xi) := \prod_{S\subseteq \{1,...,n\}}
m_S(\xi_S)$$
where $m_S\in \M(\R^{card(S)})$, the vector $\xi_S\in \R^{card(S)}$ is defined by
$\xi_S:= (\xi_i)_{i\in S}$, while $\xi \in \R^n$ is the vector $\xi:= (\xi_i)_{i=1}^n$.
Every such a symbol $m\in \M_{flag}(\R^n)$ defines naturally a generic flag paraproduct
$T_m$ by the same formula (\ref{2}). Of course, as usual, the goal is to prove
H\"{o}lder type estimates for them
\footnote{A ``flag'' is an increasing sequence of subspaces of a vector space $V$:
$\{0\}=V_0\subseteq V_1\subseteq ... \subseteq V_k = V$. It is easy to see that a generic
symbol in $\M_{flag}(\R^n)$ is singular along every possible flag of subspaces, spanned 
by the coordinate system of $\R^n$. It is also interesting to note that in a completely
different direction (see \cite{nrs}, \cite{ns}) singular integrals generated by
``flag kernels'' (this time) appear also naturally in the theory of several complex variables.}.

Let us assume, as in the case of classical paraproducts briefly discussed before, that the kernels
$K(\alpha)$, $K(\beta)$, $K(\gamma)$ are given by

$$K(\alpha) =
\sum_{k_1}
\Phi_{k_1}(\alpha_1)
\Phi_{k_1}(\alpha_2)
\Phi_{k_1}(\alpha_3),
$$

$$K(\beta) =
\sum_{k_2}
\Phi_{k_2}(\beta_1)
\Phi_{k_2}(\beta_2)
$$
and

$$K(\gamma) =
\sum_{k_3}
\Phi_{k_3}(\gamma_1)
\Phi_{k_3}(\gamma_2).
$$
In particular, the left hand side of (\ref{6}) becomes

\begin{equation}
T(f, g, h)(x) =
\sum_{k_1, k_2, k_3}
(f\ast\Phi_{k_1}\ast\Phi_{k_2})(x)\cdot
(g\ast\Phi_{k_1}\ast\Phi_{k_2}\ast\Phi_{k_3})(x)\cdot
(h\ast\Phi_{k_1}\ast\Phi_{k_3})(x)
\end{equation}
and it should be clear by looking at this expression, that there are no ``easy Banach spaces estimates''
this time. Moreover, assuming that such estimates existed, using the Calder\'{o}n - Zygmund
decomposition as before to get the ``quasi - Banach estimates'' would not help either, because
of the multi-parameter structure of the kernel $K(\alpha)K(\beta)K(\gamma)$.

In other words, completely new ideas are necessary to understand the boundedness properties of 
these flag paraproducts. More on this later on, in the last section of the paper.
We end the current one with the following result from \cite{c}.

\begin{theorem}\label{main}
Let $a, b \in \M(\R^2)$. Then, the 3 - linear operator $T_{ab}$ defined by the formula

$$T_{ab}(f_1, f_2, f_3)(x) :=
\int_{\R^3}
a(\xi_1, \xi_2)\cdot b(\xi_2, \xi_3)
\widehat{f_1}(\xi_1)
\widehat{f_2}(\xi_2)
\widehat{f_3}(\xi_3)
e^{2 \pi i x (\xi_1 + \xi_2 + \xi_3)}
d \xi$$
maps $L^{p_1}\times L^{p_2}\times L^{p_3} \rightarrow L^p$ boundedly, as long as
$1<p_1, p_2, p_3 <\infty$ and $1/p_1 + 1/p_2 + 1/p_3 =1/p$.
\end{theorem}
In addition, it has also been proven in \cite{c} that $T_{ab}$ maps also $L^{\infty}\times L^p\times L^q\rightarrow L^r$, 
$L^p\times L^{\infty}\times L^q\rightarrow L^r$, $L^p\times L^q\times L^{\infty}\rightarrow L^r$ and $L^{\infty}\times L^s\times L^{\infty}\rightarrow L^s$
boundedly, as long as $1<p,q,s<\infty$, $0<r<\infty$ and $1/p+1/q=1/r$. The only $L^{\infty}$ estimates that are not available, are those of the form
$L^{\infty}\times L^{\infty}\times L^{\infty}\rightarrow L^{\infty}$, $L^{\infty}\times L^{\infty}\times L^s\rightarrow L^s$ and
$L^s\times L^{\infty}\times L^{\infty}\rightarrow L^s$. But this should be not surprising since such estimates are in general false, as one can easily see by taking
$f_2$ to be identically equal to $1$ in the formula above.

This operator $T_{ab}$ is the simplest flag paraproduct whose complexity goes beyond the one of a Coifman - Meyer paraproduct.
However, as we remarked in \cite{c}, we believe that a similar result holds for generic
flag paraproducts of arbitrary complexity, and we plan to address this general case in a future paper \cite{c1}.

In the next three sections we will try to answer (at least partially) the question
``Why is it worthwhile to consider and study this new class of operators ?'' by
describing three distinct instances from the theory of differential equations, where they 
appear naturally.

\section{AKNS systems}

Let $\lambda\in\R$, $\lambda\neq 0$ and consider the system of differential equations

\begin{equation}\label{7}
u' = i \lambda D u + N u
\end{equation}
where $u = [u_1,...,u_n]^t$ is a vector valued function defined on the real line,
$D$ is a diagonal $n\times n$ constant matrix with real and distinct entries $d_1,...,d_n$
and $N = (a_{ij})_{i,j=1}^n$ is a matrix valued function defined also on the real line and
having the property that $a_{ii}\equiv 0$ for every $i=1,...,n$. These systems play a 
fundamental role in mathematical physics and scattering theory and they are called AKNS systems \cite{ablowitzsegur}.  The particular case
$n=2$ is also known to be deeply connected to the classical theory of Schr\"{o}dinger
operators \cite{ck1}, \cite{ck2}.

If $N\equiv 0$ it is easy to see that our system (\ref{7}) becomes a union of independent
single equations

$$u'_k = i\lambda d_k u_k$$
for $k=1,...,n$ whose solutions are

$$u_k^{\lambda}(x) = C_{k,\lambda} e^{ i\lambda d_k x}$$
and they are all $L^{\infty}(\R)$-functions. An important problem in the field is the 
following.

\begin{problem}
Prove (or disprove) that as long as $N$ is a matrix whose entries are $L^2(\R)$ functions,
then for almost every real $\lambda$, the corresponding solutions $(u_k^{\lambda})_{k=1}^n$
are all bounded functions.
\footnote{The conjecture is easy for $L^1(\R)$ entries, holds true for $L^p(\R)$ entries when
$1\leq p <2$, thanks
to the work of Christ and Kiselev \cite{ck1}, \cite{ck2} and is false for $p>2$, \cite{simon}.  }
\end{problem}

When $N\nequiv 0$ one can use a simple variation of constants argument and write $u_k(x)$
as

$$u_k(x) := e^{ i\lambda d_k x} v_k(x)$$
for $k=1,...,n$. As a consequence, the column vector $v=[v_1,...,v_n]^t$ becomes the solution
of the following system

\begin{equation}\label{8}
v' = W v
\end{equation}
where the entries of $W$ are given by $w_{lm}(x):= a_{lm}(x)e^{ i\lambda (d_l-d_m) x}$.
It is therefore enough to prove that the solutions of (\ref{8}) are bounded as long as the
entries $a_{lm}$ are square integrable.

In the particular case when the matrix $N$ is upper (or lower) triangular, the system (\ref{8}) can be solved explicitly. A straightforward
calculation shows that every single entry of the vector $v(x)$ can be written as a finite sum of expressions of the form

\begin{equation}\label{9}
\int_{x_1<...<x_k<x}
f_1(x_1) ... f_k(x_k) e^{i\lambda (\#_1 x_1 + ... + \#_k x_k)} 
d x, 
\end{equation}
where $f_1, ..., f_k$ are among the entries of the matrix $N$, while $\#_1, ..., \#_k$ are various differences of type $d_l - d_m$ as before
and satisfying the nondegeneracy condition

$$\sum_{j=j_1}^{j_2} \#_j \neq 0$$
for every $1\leq j_1 < j_2\leq k$.

Given the fact that all the entries of the matrix $N$ are $L^2(\R)$ functions and using Plancherel, one can clearly conclude that the expression
(\ref{9}) is bounded for a. e. $\lambda$, once one proves the following inequality

\begin{equation}
\left\|
\sup_M
\left|
\int_{\xi_1< ... <\xi_k<M}
\widehat{f_1}(\xi_1) ...
\widehat{f_k}(\xi_k)
e^{2\pi i x(\#_1 \xi_1 + ... +\#_k \xi_k)} 
d \xi
\right|
\right\|_{L^{2/k}_x}
\lesssim
\|f_1\|_2\cdot ... \cdot
\|f_k\|_2 < \infty.
\end{equation}
A simpler, non-maximal variant of it would be

\begin{equation}
\left\|
\int_{\xi_1< ... <\xi_k}
\widehat{f_1}(\xi_1) ...
\widehat{f_k}(\xi_k)
e^{2\pi i x(\#_1 \xi_1 + ... +\#_k \xi_k)} 
d \xi
\right\|_{L^{2/k}_x}
\lesssim
\|f_1\|_2\cdot ... \cdot
\|f_k\|_2 < \infty.
\end{equation}

The expression under the quasi-norm can be seen as a $k$-linear multiplier with symbol
$\chi_{\xi_1< ... <\xi_k} = \chi_{\xi_1<\xi_2}\cdot ... \cdot \chi_{\xi_{k-1}<\xi_k}$. Now, the ``lacunary variant'' of this multi-linear operator 
\footnote{It is customary to do this, when one faces operators which have some type of modulation invariance. For instance, the ``lacunary version'' of the Carleson operator
is the maximal Hibert transform, while the ``lacunary version'' of the bi-linear Hilbert transform is a paraproduct \cite{laceyt1}, \cite{laceyt2}. The surprise we had in
\cite{mtt:walshbiest}, \cite{mtt:fourierbiest} whith these operators is that even their ``lacunary versions'' hadn't been considered before. The reader more insterested
in learning about these operators is refered to the recent paper \cite{mtt:multiest}.}
is obtained
by replacing every bi-linear Hilbert transform type symbol $\chi_{\xi_{j-1}<\xi_j}$ \cite{laceyt1} with a smoother one $m(\xi_{j-1}, \xi_j)$, in the class $\M(\R^2)$.
The new resulted expression is clearly a flag paraproduct.

\section{General Leibnitz rules}

The following inequality of Kato and Ponce \cite{kp} plays an important role in non-linear PDEs \footnote{ $\widehat{D^{\alpha} f}(\xi) := |\xi|^{\alpha}$ for 
any $\alpha > 0$}

\begin{equation}\label{10}
\|D^{\alpha}(fg)\|_p \lesssim
\|D^{\alpha} f\|_{p_1}\|g\|_{q_1} +
\|f\|_{p_2}\|D^{\alpha} g\|_{q_2}
\end{equation}
for any $1<p_i, q_i\leq\infty$, $1/p_i + 1/q_i = 1/p$ for $i=1,2$ and $0<p<\infty$. 

It is known that the inequality holds for an arbitrary number of factors and it is also known that it is in general false if one of the indices
$p_i, q_i$ is strictly smaller than one. Given (\ref{10}), it is natural to ask if one has similar estimates for more complex expressions, such as

\begin{equation}\label{11}
\left\|
D^{\alpha}[ D^{\beta}(f_1 f_2 f_3)\cdot D^{\gamma}(f_4 f_5)]\right\|_p.
\end{equation}
Clearly, one can first apply (\ref{10}) for two factors and majorize (\ref{11}) by

\begin{equation}\label{12}
\|D^{\alpha+\beta}(f_1 f_2 f_3)\|_{p_1}\|D^{\gamma}(f_4 f_5)\|_{q_1} +
\|D^{\beta}(f_1 f_2 f_3)\|_{p_2}\|D^{\alpha+\gamma}(f_4 f_5)\|_{q_2}
\end{equation}
and after that, one can apply (\ref{10}) four more times for two and three factors, to obtain a final upper bound. However, this iterative procedure has a problem.
It doesn't work if one would like to end up with products of terms involving (for instance) only $L^2$ norms, since then one
has to have $p_1=p_2=2/3$ and $q_1=q_2 = 1$ in (\ref{12}), for which the corresponding (\ref{10}) doesn't hold.

The usual way to prove such ``paradifferential Leibnitz rules'' as the one in (\ref{10}), is by reducing them to the Coifman - Meyer theorem mentioned before.
Very briefly, the argument works as follows. First, one uses a standard Littlewood - Paley decomposition \cite{stein} and writes both $f$ and $g$ as

$$f = \sum_{k\in\Z} f\ast\Psi_k$$
and

$$g = \sum_{k\in\Z} g\ast\Psi_k$$
where $(\Psi_k)_k$ is a well chosen family of ``$\Psi$ type''. In particular, one has

$$fg = \sum_{k_1, k_2}
(f\ast\Psi_{k_1})(g\ast\Psi_{k_2}) = \sum_{k_1\sim k_2} + \sum_{k_1<<k_2} + \sum_{k_2<<k_1} :=$$

$$I + II + III.$$
Then, one can write term II (for instance) as

$$\sum_{k_1<<k_2}(f\ast\Psi_{k_1})(g\ast\Psi_{k_2}) =
\sum_{k_2}\left(\sum_{k_1<< k_2} (f\ast\Psi_{k_1})\right) (g\ast\Psi_{k_2}) :=$$

$$\sum_{k_2}(f\ast\Phi_{k_2})(g\ast\Psi_{k_2}) = \sum_{k}(f\ast\Phi_{k})(g\ast\Psi_{k}) =$$

$$\sum_k [(f\ast\Phi_k)(g\ast\Psi_k)]\ast\widetilde{\Psi}_k$$
for a well chosen family $(\widetilde{\Psi}_k)_k$ of ``$\Psi$ type'', where $(\Phi_k)_k$ is a ``$\Phi$ type'' family now.

Denote by

$$\Pi(f,g) = \sum_k [(f\ast\Phi_k)(g\ast\Psi_k)]\ast\widetilde{\Psi}_k.$$
Then, we have

$$D^{\alpha}(\Pi(f,g)) = \sum_k [(f\ast\Phi_k)(g\ast\Psi_k)]\ast D^{\alpha}\widetilde{\Psi}_k:=$$

$$\sum_k [(f\ast\Phi_k)(g\ast\Psi_k)]\ast 2^{k\alpha}\widetilde{\widetilde{\Psi}}_k =$$

$$\sum_k [(f\ast\Phi_k)(g\ast 2^{k\alpha}\Psi_k)]\ast \widetilde{\widetilde{\Psi}}_k :=$$

$$\sum_k [(f\ast\Phi_k)(g\ast D^{\alpha}\widetilde{\widetilde{\widetilde{\Psi}}}_k)]\ast \widetilde{\widetilde{\Psi}}_k =$$

$$\sum_k [(f\ast\Phi_k)(D^{\alpha} g\ast\widetilde{\widetilde{\widetilde{\Psi}}}_k)]\ast \widetilde{\widetilde{\Psi}}_k :=$$

$$\widetilde{\Pi}(f, D^{\alpha} g).$$
Now, it is easy to see that both $\Pi$ and $\widetilde{\Pi}$ are in fact bi-linear paraproducts whose symbols are

$$\sum_k\widehat{\Phi_k}(\xi_1)\widehat{\Psi_k}(\xi_2)\widehat{\widetilde{\Psi_k}}(\xi_1+\xi_2)$$
and

$$\sum_k\widehat{\Phi_k}(\xi_1)\widehat{\widetilde{\widetilde{\widetilde{\Psi_k}}}}(\xi_2)\widehat{\widetilde{\widetilde{\Psi_k}}}(\xi_1+\xi_2)$$
respectively.

Combining the above equality (between $\Pi$ and $\widetilde{\Pi}$) with the Coifman - Meyer theorem and treating similarly the other two terms I and III, give
the desired (\ref{10}).

What we've learned from the calculations above, is that evey non-linearity of the form $D^{\alpha}(fg)$ can be mollified and written as a finite sum of various 
bi-linear paraproducts applied to either $D^{\alpha}f$ and $g$ or to $f$ and $D^{\alpha}g$ and this fact allows one to reduce inequality (\ref{10}) to
the Coifman - Meyer theorem on paraproducts.

Let us consider now the non-linear expression $D^{\alpha}(D^{\beta}(fg) h)$ which is the simplest non-linearity of the same complexity as the one in (\ref{11}).
Clearly, this non-linearity which we say is of complexity 2, can be seen as a composition of two non-linearities of complexity 1. In particular, 
this may suggest that one way to mollify it
is by composing the mollified versions of its lower complexity counterparts, thus obtaining expressions of type
$\Pi'(\Pi''(F,G), H)$ where $\Pi'$ and $\Pi''$ are bi-linear paraproducts as before.
This procedure reduces the problem of estimating $\|D^{\alpha}(D^{\beta}(fg) h)\|_p$ to the problem of estimating $\|\Pi'(\Pi''(F,G), H)\|_p$, which can be done
by applying the Coifman - Meyer theorem two times in a row. However, as we pointed out before, this point of view cannot handle the case when for instance both
$F$ and $G$ are functions in an $L^{1+\epsilon}$ space for $\epsilon > 0$ a small number. To be able to understand completely these non-linearities of higher 
complexity one has to proceed differently. As we will see, the correct point of view is not to write $D^{\alpha}(D^{\beta}(fg) h)$ as a sum of composition of
paraproducts, but instead as a sum of flag paraproducts.

Clearly, in order to achieve this, we need to understand two things. First, how to mollify an expression of type $\Pi(F,G)H$ and then, how to mollify
$D^{\alpha} (\Pi(F,G)H)$.

Let us assume as before that $\Pi(F,G)$ is given by

$$\Pi(F,G) = \sum_{k_1<<k_2} (F\ast\Psi_{k_1}) (G\ast\Psi_{k_2})$$
and decompose $H$ as usual as

$$H = \sum_{k_3} H\ast\Psi_{k_3}.$$
As a consequence, we have

$$\Pi(F,G) H =$$

$$\sum_{k_1<<k_2; k_3} (F\ast\Psi_{k_1}) (G\ast\Psi_{k_2}) (H\ast\Psi_{k_3}) =$$

$$\sum_{k_1<<k_2; k_3<< k_2} + \sum_{k_1<<k_2; k_3 \sim k_2} + \sum_{k_1<<k_2<<k_3}:= $$

$$A + B + C.$$
It is not difficult to see that both $A$ and $B$ are simply 3-linear paraproducts, while $C$ can be written as

$$\int_{\R^3}
m(\xi_1, \xi_2, \xi_3)
\widehat{F}(\xi_1)
\widehat{G}(\xi_2)
\widehat{H}(\xi_3)
e^{2\pi i x(\xi_1 +\xi_2 +\xi_3)} 
d \xi
$$
where

$$m(\xi_1, \xi_2, \xi_3) = \sum_{k_1<<k_2<<k_3}\widehat{\Psi_{k_1}}(\xi_1)\widehat{\Psi_{k_2}}(\xi_2)\widehat{\Psi_{k_3}}(\xi_3) :=$$

$$\sum_{k_2<<k_3}\widehat{\Phi_{k_2}}(\xi_1)\widehat{\Psi_{k_2}}(\xi_2)\widehat{\Psi_{k_3}}(\xi_3) = $$

\begin{equation}\label{13}
\sum_{k_2<<k_3}\widehat{\Phi_{k_2}}(\xi_1)\widehat{\Psi_{k_2}}(\xi_2) \widehat{\widetilde{\Phi_{k_3}}}(\xi_2) \widehat{\Psi_{k_3}}(\xi_3)
\end{equation}
for some well chosen family $(\widetilde{\Phi_{k_3}}(\xi_2))_{k_3}$ of ``$\Phi$ type''. But then, one observes that (\ref{13}) splits as

$$\left(
\sum_{k_2}\widehat{\Phi_{k_2}}(\xi_1)\widehat{\Psi_{k_2}}(\xi_2)\right)
\left(\sum_{k_3}\widehat{\widetilde{\Phi_{k_3}}}(\xi_2) \widehat{\Psi_{k_3}}(\xi_3)\right)
$$
since the only way in which $\widehat{\Psi_{k_2}}(\xi_2) \widehat{\widetilde{\Phi_{k_3}}}(\xi_2)\neq 0$ is to have $k_2<<k_3$.
This shows that the symbol of $C$ can be written as $a(\xi_1, \xi_2) b(\xi_2,\xi_3)$ with both $a$ and $b$ in $\M(\R^2)$, which means that $C$ is indeed a flag paraproduct.

In order to understand now how to molify $D^{\alpha} (\Pi(F,G)H)$, let us first rewrite (\ref{13}) as

\begin{equation}\label{14}
\sum_{k_2<<k_3}\widehat{\Phi_{k_2}}(\xi_1)\widehat{\Psi_{k_2}}(\xi_2) \widehat{\widetilde{\Phi_{k_2}}}(\xi_1+\xi_2) 
\widehat{\widetilde{\widetilde{\Phi_{k_3}}}}(\xi_1+\xi_2) 
\widehat{\Psi_{k_3}}(\xi_3)
\widehat{\widetilde{\Phi_{k_3}}}(\xi_1+\xi_2+\xi_3)
\end{equation}
where as usual, $\widehat{\widetilde{\Phi_{k_2}}}(\xi_1+\xi_2)$, $\widehat{\widetilde{\widetilde{\Phi_{k_3}}}}(\xi_1+\xi_2)$ and   
$\widehat{\widetilde{\Phi_{k_3}}}(\xi_1+\xi_2+\xi_3)$ have been inserted naturally into (\ref{13}).

The advantage of (\ref{14}) is that the corresponding 3-linear operator can be easily written as

\begin{equation}\label{15}
\sum_{k_3}
\left\{
\left(
\sum_{k_2<<k_3}
[(F\ast\Phi_{k_2})(G\ast\Psi_{k_2})]\ast\widetilde{\Phi_{k_2}}\right)\ast\widetilde{\widetilde{\Phi_{k_3}}} \cdot (F\ast\Psi_{k_3})\right\}\ast
\widetilde{\Phi_{k_3}}.
\end{equation}
Then, exactly as before one can write

$$D^{\alpha} (\Pi(F,G)H) =$$

$$\sum_{k_3}
\left\{
\left(
\sum_{k_2<<k_3}
[(F\ast\Phi_{k_2})(G\ast\Psi_{k_2})]\ast\widetilde{\Phi_{k_2}}\right)\ast\widetilde{\widetilde{\Phi_{k_3}}} \cdot (F\ast\Psi_{k_3})\right\}\ast
D^{\alpha}\widetilde{\Phi_{k_3}} =$$

$$= ... =$$

$$
\sum_{k_3}
\left\{
\left(
\sum_{k_2<<k_3}
[(F\ast\Phi_{k_2})(G\ast\Psi_{k_2})]\ast\widetilde{\Phi_{k_2}}\right)\ast\widetilde{\widetilde{\Phi_{k_3}}} 
\cdot (D^{\alpha}F\ast\widetilde{\Psi_{k_3}})\right\}\ast
\widetilde{\widetilde{\Phi_{k_3}}}.
$$
Now, this tri-linear operator has the symbol

$$\sum_{k_2<<k_3}
\widehat{\Phi_{k_2}}(\xi_1)
\widehat{\Psi_{k_2}}(\xi_2)
\widehat{\widetilde{\Psi_{k_2}}}(\xi_1+\xi_2)
\widehat{\widetilde{\widetilde{\Psi_{k_3}}}}(\xi_1+\xi_2)
\widehat{\widetilde{\Psi_{k_3}}}(\xi_3)
\widehat{\widetilde{\widetilde{\Psi_{k_3}}}}(\xi_1+\xi_2 +\xi_3) =
$$

$$
\sum_{k_2<<k_3}
\widehat{\Phi_{k_2}}(\xi_1)
\widehat{\Psi_{k_2}}(\xi_2)
\widehat{\widetilde{\Psi_{k_3}}}(\xi_3)
\widehat{\widetilde{\widetilde{\Psi_{k_3}}}}(\xi_1+\xi_2 +\xi_3) =
$$

$$
\sum_{k_2<<k_3}
\widehat{\Phi_{k_2}}(\xi_1)
\widehat{\Psi_{k_2}}(\xi_2)
\widehat{\widetilde{\Phi_{k_3}}}(\xi_2)
\widehat{\widetilde{\Psi_{k_3}}}(\xi_3)
\widehat{\widetilde{\widetilde{\Psi_{k_3}}}}(\xi_1+\xi_2 +\xi_3) 
$$
and this splits again as

$$
\left(\sum_{k_2}
\widehat{\Phi_{k_2}}(\xi_1)
\widehat{\Psi_{k_2}}(\xi_2)\right)
\left(\sum_{k_3}
\widehat{\widetilde{\Phi_{k_3}}}(\xi_2)
\widehat{\widetilde{\Psi_{k_3}}}(\xi_3)
\widehat{\widetilde{\widetilde{\Psi_{k_3}}}}(\xi_1+\xi_2 +\xi_3)\right) :=
$$

$$\alpha(\xi_1,\xi_2)\cdot \beta(\xi_1, \xi_2, \xi_3)$$
which is an element of $\M_{flag}(\R^3)$.

Since tri-linear operators of the form (\ref{15}) give rise to model operators which have been understood in \cite{c}, the above discussion proves the following
``grand Leibnitz rule'' for the simplest non-linearity of complexity 2

$$\|D^{\alpha}(D^{\beta}(fg) h)\|_p \lesssim$$

$$
\|D^{\alpha+\beta}f\|_{p_1} \|g\|_{q_1} \|h\|_{r_1} +
\|f\|_{p_2} \|D^{\alpha+\beta}g\|_{q_2} \|h\|_{r_2} +
\|D^{\beta}f\|_{p_3}\|g\|_{q_2}\|D^{\alpha}h\|_{r_3} +
\|f\|_{p_4}\|D^{\beta}g\|_{q_4}\|D^{\alpha}h\|_{r_4}$$
for every $1<p_i, q_i, r_i<\infty$ with $1/p_i+1/q_i+1/r_i = 1/p$ for $i=1,2,3,4$.

\section{Flag Paraproducts and the non-linear Schr\"{o}dinger equation}

In this section the goal is to briefly describe some recent study of Germain, Masmoudi and Shatah \cite{gms}
on the non-linear Schr\"{o}dinger equation. We are grateful to Pierre Germain who explained
parts of this work to us. 

In just a few words, the main task of these three authors is to develop a general method for
understanding the global existence for small initial data of various non-linear 
Schr\"{o}dinger equations and systems of non-linear Schr\"{o}dinger equations.

Consider the following ``quadratic example''

\begin{equation}
\partial_t u + i\triangle u = u^2
\end{equation}

\begin{equation}
u|_{t=2} = u_2
\end{equation}
for $(t,x)\in \R\times \R^n$
\footnote{There are some technical reasons for which the authors prefered the initial time
to be 2, related to the norm of the space $X$ where the global existence takes place.}.

The Duhamel formula written on the Fourier side becomes

$$\widehat{u} (t,\xi) =
\widehat{u_2}(\xi) e^{i t |\xi|^2} + \int_2^t
e^{-i (s-t)|\xi|^2} \widehat{u^2}(s,\xi) d s.
$$

Then, if one writes $u=e^{-it\triangle} f$, the above formula becomes

\begin{equation}\label{16}
\widehat{f}(t,\xi) =
\widehat{u_2}(\xi) + \int_2^t \int
e^{is (-|\xi|^2 + |\eta|^2 + |\xi-\eta|^2)} 
\widehat{f}(s,\eta)\widehat{f}(s,\xi-\eta)
d s
d \eta.
\end{equation}

A significant part of the argument in \cite{gms} depends on how well one estimates
the integral expression in (\ref{16}). The idea would be to take advantage of the oscillation
of the term $e^{is\Phi}$ where $\Phi:= -|\xi|^2 + |\eta|^2 + |\xi-\eta|^2$.
However, $\Phi$ is ``too degenerate'' (i.e. $\Phi=0$ whenever $\xi=\eta$ or $\eta=0$)
and as a consequence, the usual ``$\frac{d}{ds}\{\frac{e^{is\Phi}}{i\Phi}\}$'' argument
doesn't work. One needs a ``wiser integration by parts'', not only in the $s$ variable but also
in the $\eta$ variable. Denote by

$$P:= -\eta +\frac{1}{2}\xi$$
and

$$Z:= \Phi + P\cdot (\partial_{\eta}\Phi).
$$
Observe that 

$$Z = - (|\eta|^2 + |\xi-\eta|^2)
$$
which is identically equal to zero only when $\xi=\eta=0$.

Alternatively, one has

$$\frac{1}{iZ} (\partial_s + \frac{P}{s}\partial_{\eta}) e^{is\Phi} = e^{is\Phi}.$$
In particular, the inverse Fourier transform of the integral term in (\ref{16}) becomes

$$\F^{-1} \int_2^t\int
\frac{1}{iZ} (\partial_s + \frac{P}{s}\partial_{\eta}) e^{is\Phi}
\widehat{f}(s,\xi-\eta)
\widehat{f}(s,\eta)
d s
d \eta =
$$

$$I + II.$$
Using the fact that $|\xi-\eta|^2/iZ$ and $|\eta|^2/iZ$ are both classical symbols in
$\M(\R^{2n})$, Coifman - Meyer theorem proves that both $I$ and $II$ are ``smoothing expressions''.
To conclude, expressions of type

$$\F^{-1}\int_2^t\int
e^{is\Phi}
m(\xi,\eta)
\widehat{g}(s,\eta)
\widehat{h}(s,\xi-\eta)
d s
d \eta
$$
for some $m\in \M(\R^{2n})$ appear naturally, and it is easy to see that if one keeps $s$ fixed, the rest of the formula
is just a bi-linear paraproduct.

Coming back to $I$, one of the expressions related to it (after the integration by parts)
is of the form

\begin{equation}\label{17}
\F^{-1}\int_2^t\int
e^{is\Phi}
m(\xi,\eta)
\partial_s\widehat{f}(s,\eta)
\widehat{F}(s,\xi-\eta)
d s
d \eta
\end{equation}
for a certain new function $F$.
Since $u=e^{-is\triangle} f$ it follows that $f = e^{is\triangle} u$ and then, by using the 
fact that $u$ solves the equation, we deduce that $\partial_s f = e^{is\triangle} u^2$
which means that

$$\widehat{\partial_s f}(s,\eta) = e^{-is|\eta|^2}\widehat{u^2}(\eta) =
$$

$$e^{-is|\eta|^2}
\int_{\eta_1+\eta_2 =\eta}
\widehat{u}(\eta_1)
\widehat{u}(\eta_2)
d \eta_1
d \eta_2 =$$

$$e^{-is|\eta|^2}
\int_{\eta_1+\eta_2 =\eta}
e^{is|\eta_1|^2}
\widehat{f}(s,\eta_1)
e^{is|\eta_2|^2}
\widehat{f}(s,\eta_2)
d \eta_1
d \eta_2 =$$

$$e^{-is|\eta|^2}
\int
e^{is|\tau|^2}
\widehat{f}(s,\tau)
e^{is|\eta-\tau|^2}
\widehat{f}(s,\eta-\tau)
d \tau.$$
Using this in (\ref{17}) one obtains an expression of the form

$$
\F^{-1}
\int_2^t\int
e^{is\widetilde{\Phi}}
m(\xi, \eta)
\widehat{f}(s,\tau)
\widehat{f}(s, \eta-\tau)
\widehat{F}(s,\xi-\eta)
d s
d \eta
d \tau
$$
where $\widetilde{\Phi}:= -|\xi|^2 + |\xi-\eta|^2 + |\tau|^2 + |\eta-\tau|^2$.

Using now a similar ``integration by parts argument'' in all three variables this time
one obtains as before expressions of type

$$
\F^{-1}
\int_2^t\int
e^{is\widetilde{\Phi}}
m(\xi, \eta) m(\xi,\eta,\tau)
\widehat{f}(s,\tau)
\widehat{f}(s, \eta-\tau)
\widehat{F}(s,\xi-\eta)
d s
d \eta
d \tau.
$$
The inner formula (for a fixed $s$) is of the form

$$\F^{-1}
\int
m(\xi, \eta) m(\xi,\eta,\tau)
\widehat{g}(\tau)
\widehat{f}(\eta-\tau)
\widehat{h}(\xi-\eta)
d s
d \eta
d \tau
$$
and we claim that it is naturally related to the flag paraproducts we described earlier.

Alternatively, one can rewrite it as

\begin{equation}\label{18}
\int_{\R^3}
m(\xi, \eta) m(\xi,\eta,\gamma)
\widehat{f}(\eta-\gamma)
\widehat{g}(\gamma)
\widehat{h}(\xi-\eta)
e^{ix\xi}
d \xi
d \eta
d \gamma.
\end{equation}
Then, if we change variables $\xi-\eta:=\xi_3$, $\eta-\gamma:=\xi_1$ and $\gamma:=\xi_2$
(\ref{18}) becomes

$$\int_{\R^3}
m(\xi_1+\xi_2+\xi_3,\xi_1+\xi_2) m(\xi_1+\xi_2+\xi_2, \xi_1+\xi_2,\xi_2)
\widehat{f}(\xi_1)
\widehat{g}(\xi_2)
\widehat{h}(\xi_3)
e^{i x (\xi_1+\xi_2+\xi_3)}
d \xi :=
$$

\begin{equation}\label{19}
\int_{\R^3}
\widetilde{m}(\xi_1+\xi_2,\xi_3) \widetilde{\widetilde{m}}(\xi_1, \xi_2, \xi_3)
\widehat{f}(\xi_1)
\widehat{g}(\xi_2)
\widehat{h}(\xi_3)
e^{i x (\xi_1+\xi_2+\xi_3)}
d \xi
\end{equation}
where $\widetilde{m}\in\M(\R^2)$ while $\widetilde{\widetilde{m}}\in\M(\R^3)$.

As it stands, (\ref{19}) is not a flag paraproduct, but we will show that its analysis
can be reduced to the analysis of a flag paraproduct.

Assume for simplicity that we have

$$
\widetilde{m}(\xi_1+\xi_2,\xi_3) = \sum_{k_1} 
\widehat{\Phi_{k_1}}(\xi_1+\xi_2)
\widehat{\Phi_{k_1}}(\xi_3)
$$
and

$$
\widetilde{\widetilde{m}}(\xi_1, \xi_2, \xi_3) =
\sum_{k_2}
\widehat{\Phi_{k_2}}(\xi_1)
\widehat{\Phi_{k_2}}(\xi_2)
\widehat{\Phi_{k_2}}(\xi_3)
$$
as before\footnote{Of course, everything is defined in $\R^n$ now, but the extensions to arbitrary euclidean spaces are straightforward}. Then,

$$
\widetilde{m}(\xi_1+\xi_2,\xi_3)
\widetilde{\widetilde{m}}(\xi_1, \xi_2, \xi_3) =
$$

\begin{equation}\label{20}
\sum_{k_1, k_2}
\widehat{\Phi_{k_1}}(\xi_1+\xi_2)
\widehat{\Phi_{k_1}}(\xi_3)
\widehat{\Phi_{k_2}}(\xi_1)
\widehat{\Phi_{k_2}}(\xi_2)
\widehat{\Phi_{k_2}}(\xi_3).
\end{equation}
Clearly, we have two interesting cases:

\underline{Case 1: $k_2<<k_1$}.

Here, the only possibility is to have $(\widehat{\Phi_{k_1}}(\xi_3))_{k_1}$
of $\Phi$ type in which case $(\widehat{\Phi_{k_1}}(\xi_1+\xi_2))_{k_1}$ must be
of ``$\Psi$ type''. Since (\ref{20}) can also be written as

$$
\sum_{k_2<< k_1}
\widehat{\Phi_{k_1}}(\xi_1+\xi_2)
\widehat{\Phi_{k_1}}(\xi_3)
\widehat{\Phi_{k_2}}(\xi_1)
\widehat{\Phi_{k_2}}(\xi_2)
\widehat{\widetilde{\Phi_{k_2}}}(\xi_1+\xi_2)
\widehat{\Phi_{k_2}}(\xi_3)
$$
for a well chosen family $(\widehat{\widetilde{\Phi_{k_2}}}(\xi_1+\xi_2))_{k_2}$
we see that the only way in which 
$\widehat{\Phi_{k_1}}(\xi_1+\xi_2)
\widehat{\widetilde{\Phi_{k_2}}}(\xi_1+\xi_2) \neq 0$ is to have $k_1\sim k_2$. But in this 
case, the multiplier belongs to $\M(\R^3)$ and we simply face a tri-linear paraproduct.

\underline{Case 2: $k_2>>k_1$}.

This time, the only possibility is to have $(\widehat{\Phi_{k_2}}(\xi_3))_{k_2}$ of
``$\Phi$'' type. Then, we can ``complete''the expression in (\ref{20}) as

\begin{equation}\label{21}
\sum_{k_2>>k_1}
\widehat{\Phi_{k_2}}(\xi_1)
\widehat{\Phi_{k_2}}(\xi_2)
\widehat{\Phi_{k_2}}(\xi_3)
\widehat{\widetilde{\Phi_{k_2}}}(\xi_1+\xi_2)
\widehat{\Phi_{k_1}}(\xi_1+\xi_2)
\widehat{\Phi_{k_1}}(\xi_3)
\widehat{\widetilde{\Phi_{k_1}}}(\xi_1+\xi_2+\xi_3).
\end{equation}

Now, for $(\widehat{\widetilde{\Phi_{k_2}}}(\xi_1+\xi_2))_{k_2}$ we have two options. Either it is of ``$\Psi$ type'', in which situation the only non-zero case would be
when $k_1\sim k_2$. But then, this means that we are again in a paraproduct setting. Or,  $(\widehat{\widetilde{\Phi_{k_2}}}(\xi_1+\xi_2))_{k_2}$
is of ``$\Phi$ type'' but this can only happen when both $(\widehat{\Phi_{k_2}}(\xi_1))_{k_2}$ and $(\widehat{\Phi_{k_2}}(\xi_2))_{k_2}$ are of ``$\Psi$ type''
(and their oscillations cancel out). Since we also know that either $(\widehat{\Phi_{k_1}}(\xi_1+\xi_2))_{k_1}$ or $(\widehat{\Phi_{k_1}}(\xi_3))_{k_1}$
has to be of a ``$\Psi$ type'', it follows that (\ref{21}) splits as

$$\left(\sum_{k_2}\widehat{\Phi_{k_2}}(\xi_1)
\widehat{\Phi_{k_2}}(\xi_2)
\widehat{\Phi_{k_2}}(\xi_3)
\widehat{\widetilde{\Phi_{k_2}}}(\xi_1+\xi_2)\right)
\left(\sum_{k_1}
\widehat{\Phi_{k_1}}(\xi_1+\xi_2)
\widehat{\Phi_{k_1}}(\xi_3)
\widehat{\widetilde{\Phi_{k_1}}}(\xi_1+\xi_2+\xi_3)
\right).
$$
But then, if we denote by $T(f,g,h)$ the corresponding tri-linear operator, one has

$$\Lambda_T(f,g,h,k):= \int T(f,g,h)(x) k(x) d x =$$

$$\int
\left(\sum_{k_2}\widehat{\Phi_{k_2}}(\xi_1)
\widehat{\Phi_{k_2}}(\xi_2)
\widehat{\Phi_{k_2}}(\xi_3)
\widehat{\widetilde{\Phi_{k_2}}}(\xi_1+\xi_2)\right)
\left(\sum_{k_1}
\widehat{\Phi_{k_1}}(\xi_1+\xi_2)
\widehat{\Phi_{k_1}}(\xi_3)
\widehat{\widetilde{\Phi_{k_1}}}(\xi_1+\xi_2+\xi_3)
\right)\cdot
$$

$$\cdot\widehat{f}(\xi_1)
\widehat{g}(\xi_2)
\widehat{h}(\xi_3)
\widehat{k}(-\xi_1-\xi_2-\xi_3)
d \xi =
$$

$$\int
\left(\sum_{k_2}\widehat{\widetilde{\Phi_{k_2}}}(-\xi_1)
\widehat{\Phi_{k_2}}(\xi_2)
\widehat{\Phi_{k_2}}(\xi_3)
\widehat{\widetilde{\widetilde{\Phi_{k_2}}}}(-\xi_1-\xi_2)\right)
\left(\sum_{k_1}
\widehat{\widetilde{\Phi_{k_1}}}(-\xi_1-\xi_2)
\widehat{\Phi_{k_1}}(\xi_3)
\widehat{\widetilde{\widetilde{\Phi_{k_1}}}}(-\xi_1-\xi_2-\xi_3)
\right)\cdot
$$

$$\cdot\widehat{f}(\xi_1)
\widehat{g}(\xi_2)
\widehat{h}(\xi_3)
\widehat{k}(-\xi_1-\xi_2-\xi_3)
d \xi.
$$
Then, if we denote by $\lambda:= -\xi_1-\xi_2-\xi_3$ the previous expression becomes

$$\int
\left(\sum_{k_2}\widehat{\widetilde{\Phi_{k_2}}}(\lambda+\xi_2+\xi_3)
\widehat{\Phi_{k_2}}(\xi_2)
\widehat{\Phi_{k_2}}(\xi_3)
\widehat{\widetilde{\widetilde{\Phi_{k_2}}}}(\lambda+\xi_3)\right)
\left(\sum_{k_1}
\widehat{\widetilde{\Phi_{k_1}}}(\lambda+\xi_3)
\widehat{\Phi_{k_1}}(\xi_3)
\widehat{\widetilde{\widetilde{\Phi_{k_1}}}}(\lambda)
\right)\cdot
$$

$$\cdot\widehat{f}(-\lambda-\xi_2-\xi_3)
\widehat{g}(\xi_2)
\widehat{h}(\xi_3)
\widehat{k}(\lambda)
d \xi
d \lambda :=
$$

$$\int m(\xi_2,\xi_3, \lambda)
m(\xi_3,\lambda)
\widehat{g}(\xi_2)
\widehat{h}(\xi_3)
\widehat{k}(\lambda)
\widehat{f}(-\lambda-\xi_2-\xi_3)
d \xi
d \lambda :=$$

$$\int
\Pi_{flag}(g, h, k)(x) f(x) d x.
$$
In conclusion, there exists a flag paraproduct $\Pi_{flag}$ so that

$$\Lambda_T(f,g,h,k) =
\int \Pi_{flag}(g, h, k)(x) f(x) d x
$$
which reduces the study of $T$ to the study of $\Pi_{flag}$. 

As far as we understood from \cite{gms}, paraproducts appear in the study of the
$3D$ quadratic NLS, while the flag paraproducts are in addition necessary to deal with the 
more delicate $2D$ quadratic NLS.

\section{Remarks about the proof of Theorem \ref{main}}

Let us first recall that the symbol of the operator in question is $a(\xi_1,\xi_2) b(\xi_2,\xi_3)$. Split as before both $a$ and $b$ as

$$a(\xi_1,\xi_2) = \sum_{k_1} 
\widehat{\Phi_{k_1}}(\xi_1)
\widehat{\Phi_{k_1}}(\xi_2)
$$
and

$$b(\xi_2,\xi_3) = \sum_{k_2} 
\widehat{\Phi_{k_2}}(\xi_2)
\widehat{\Phi_{k_2}}(\xi_3)
$$
which means that

$$a(\xi_1,\xi_2) b(\xi_2,\xi_3) =
\sum_{k_1, k_2}
\widehat{\Phi_{k_1}}(\xi_1)
\widehat{\Phi_{k_1}}(\xi_2)
\widehat{\Phi_{k_2}}(\xi_2)
\widehat{\Phi_{k_2}}(\xi_3).
$$

As usual, there are three cases. Either $k_1\sim k_2$ or $k_1<<k_2$ or $k_2<<k_1$. The first one is easy, since it generates paraproducts and so we only need to
deal with the second one, the third one being completly symmetric. The corresponding symbol is

\begin{equation}\label{22}
\sum_{k_1<<k_2}
\widehat{\Phi_{k_1}}(\xi_1)
\widehat{\Phi_{k_1}}(\xi_2)
\widehat{\Phi_{k_2}}(\xi_2)
\widehat{\Phi_{k_2}}(\xi_3).
\end{equation}
Clearly, in this case we must have $(\widehat{\Phi_{k_2}}(\xi_2))_{k_2}$ of ``$\Phi$ type''. For reasons that will be clearer later on, we would have liked
instead of (\ref{22}) to face an expression of the form

\begin{equation}\label{23}
\sum_{k_1<<k_2}
\widehat{\Phi_{k_1}}(\xi_1)
\widehat{\Phi_{k_1}}(\xi_2)
\widehat{\Phi_{k_2}}(\xi_1+\xi_2)
\widehat{\Phi_{k_2}}(\xi_3).
\end{equation}

Indeed \footnote{We should emphasize here that a similar problem appeared in the ``bi-est case'' \cite{mtt:walshbiest}, \cite{mtt:fourierbiest}. 
There, the solution came from the 
observation that inside a region of the form $|\xi_1-\xi_2|<< |\xi_2-\xi_3|$, one simply has the equality 
$\chi_{\xi_1<\xi_2}\chi_{\xi_2<\xi_3} =  \chi_{\xi_1<\xi_2}\chi_{\frac{\xi_1+\xi_2}{2}<\xi_3}$. However, in our case a similar formula is not available 
unfortunately, since we are working 
with generic multipliers. This is the main reason for which we will have to face later on not only discrete models of type (\ref{24}) (which are the ``lacunary'' versions
of the model operators in  \cite{mtt:walshbiest}, \cite{mtt:fourierbiest}), but also the new ones described in (\ref{25}).    }, if we had to deal with (\ref{23}) instead, 
we would have completed it as

\begin{equation}
\sum_{k_1<<k_2}
\widehat{\Phi_{k_1}}(\xi_1)
\widehat{\Phi_{k_1}}(\xi_2)
\widehat{\widetilde{\Phi_{k_1}}}(\xi_1+\xi_2)
\widehat{\Phi_{k_2}}(\xi_1+\xi_2)
\widehat{\Phi_{k_2}}(\xi_3).
\widehat{\widetilde{\Phi_{k_2}}}(\xi_1+\xi_2+\xi_3)
\end{equation}
and then the 3-linear operator having this symbol could be coveniently rewritten as

$$\sum_{k_2}
\left\{
\left(
\sum_{k_1<<k_2}
[(f\ast\Phi_{k_1})(g\ast\Phi_{k_1})]\ast\widetilde{\Phi_{k_1}}\right)\ast\Phi_{k_2}\cdot (h\ast\Phi_{k_2})
\right\}\ast\widetilde{\Phi_{k_2}}
$$
which could be further discretized to become an average of model operators of the form

\begin{equation}\label{24}
\sum_I \frac{1}{|I|^{1/2}}
\langle B_I(f,g), \Phi^1_I\rangle
\langle h, \Phi^2_I\rangle
\Phi^3_I,
\end{equation}
where

$$B_I(f,g) = \sum_{J: |J|>|I|}
\frac{1}{|J|^{1/2}}
\langle f, \Phi^1_J\rangle
\langle g, \Phi^2_J\rangle
\Phi^3_J,
$$
while $\Phi^i_I$ and $\Phi^i_J$ are all $L^2$ - normalized bump functions, adapted to dyadic intervals $I$ and $J$ respectively,
which are either of ``$\Phi$ or $\Psi$ type'' \cite{c} \footnote{In fact, at least two of the families $(\Phi^i_I)_I$, $i=1,2,3$ and at least
two of the families $(\Phi^i_J)_J$, $i=1,2,3$ are of a ``$\Psi$ type''.}.

Our next goal is to explain how can one bring the original (\ref{22}) symbol to its better variant (\ref{23}).

The idea is that in the case when $k_1<<k_2$, then $\widehat{\Phi_{k_2}}(\xi_2)$ becomes very close to $\widehat{\Phi_{k_2}}(\xi_1+\xi_2)$,
since $\xi_1$ must be in the support of $\widehat{\Phi_{k_1}}(\xi_1)$, which lives at a much smaller scale.

Therefore, using a Taylor decomposition, one can write

$$\widehat{\Phi_{k_2}}(\xi_2) = \widehat{\Phi_{k_2}}(\xi_1+ \xi_2) +
\frac{\widehat{\Phi_{k_2}}'(\xi_1+ \xi_2)}{1!}(-\xi_1) + 
\frac{\widehat{\Phi_{k_2}}''(\xi_1+ \xi_2)}{2!}(-\xi_1)^2 + ... +
\frac{\widehat{\Phi_{k_2}}^M(\xi_1+ \xi_2)}{M!}(-\xi_1)^M +
R_M(\xi_1, \xi_2).$$
Clearly, the $0$th term gives rise to a model operator similar to the one before. Fix now $0<l\leq M$ and consider an intermediate term
$\frac{\widehat{\Phi_{k_2}}^l(\xi_1+ \xi_2)}{l!}(-\xi_1)^l$. If we insert it into (\ref{22}) the corresponding multiplier becomes

$$
\sum_{k_1<<k_2}
\widehat{\Phi_{k_1}}(\xi_1)
\widehat{\Phi_{k_1}}(\xi_2)
\frac{\widehat{\Phi_{k_2}}^l(\xi_1+ \xi_2)}{l!}(-\xi_1)^l
\widehat{\Phi_{k_2}}(\xi_3) =
$$

$$
\sum_{\# = 1000}^{\infty}
\sum_{k_2 = k_1 + \#}^{\infty}
\widehat{\Phi_{k_1}}(\xi_1)
\widehat{\Phi_{k_1}}(\xi_2)
\frac{\widehat{\Phi_{k_2}}^l(\xi_1+ \xi_2)}{l!}(-\xi_1)^l
\widehat{\Phi_{k_2}}(\xi_3) =
$$

$$
\sum_{\# = 1000}^{\infty}
\sum_{k_2 = k_1 + \#}^{\infty}
\widehat{\Phi_{k_1}}(\xi_1)(-\xi_1)^l
\widehat{\Phi_{k_1}}(\xi_2)
\frac{\widehat{\Phi_{k_2}}^l(\xi_1+ \xi_2)}{l!}
\widehat{\Phi_{k_2}}(\xi_3) :=
$$

$$
\sum_{\# = 1000}^{\infty}
\sum_{k_2 = k_1 + \#}^{\infty}
\frac{2^{k_1 l}}{2^{ k_2 l}}
\widehat{\Phi_{k_1, l}}(\xi_1)
\widehat{\Phi_{k_1}}(\xi_2)
\widehat{\Phi_{k_2, l}}(\xi_1+ \xi_2)
\widehat{\Phi_{k_2}}(\xi_3) =
$$

$$
\sum_{\# = 1000}^{\infty} 2^{-\# l}
\sum_{k_2 = k_1 + \#}^{\infty}
\widehat{\Phi_{k_1, l}}(\xi_1)
\widehat{\Phi_{k_1}}(\xi_2)
\widehat{\Phi_{k_2, l}}(\xi_1+ \xi_2)
\widehat{\Phi_{k_2}}(\xi_3).
$$
But then, the inner expression (for a fixed $\#$) generates model operators of the form (see \cite{c})

\begin{equation}\label{25}
\sum_I \frac{1}{|I|^{1/2}}
\langle B_I^{\#}(f,g), \Phi^1_I\rangle
\langle h, \Phi^2_I\rangle
\Phi^3_I
\end{equation}
where

$$B_I^{\#}(f,g) = \sum_{J: |J|= 2^{\#}|I|}
\frac{1}{|J|^{1/2}}
\langle f, \Phi^1_J\rangle
\langle g, \Phi^2_J\rangle
\Phi^3_J.
$$
Finally, it is not difficult to see that the ``rest'' operator (the one which corresponds to $R_M(\xi_1,\xi_2)$) can be written as

$$\sum_{\# = 1000}^{\infty} 2^{-\# M}
T_{\#}(f,g,h)$$
where the symbol $m_{\#}$ of the operator $T_{\#}$ satisfies an estimate of the type

$$|\partial^{\alpha}m_{\#}(\xi)|\lesssim
2^{\# |\alpha|}\frac{1}{|\xi|^{|\alpha|}}$$
for many multi-indices $\alpha$. As a consequence, Coifman - Meyer theorem implies that each $T_{\#}$ is bounded with a bound of type $O(2^{100 \#})$,
which is acceptable if we pick $M$ large enough. 

All of these show that one needs to understand the model operators (\ref{24}) and (\ref{25}), in order to prove our theorem.

The 4-linear form associated to (\ref{24}) can be written as

$$\Lambda(f,g,h,k) =
$$

$$
\sum_I \frac{1}{|I|^{1/2}}
\langle B_I(f,g), \Phi^1_I\rangle
\langle h, \Phi^2_I\rangle
\langle k, \Phi^3_I\rangle =
$$

$$\sum_J
\frac{1}{|J|^{1/2}}
\langle f, \Phi^1_J\rangle
\langle g,\Phi^2_J\rangle
\langle B_J(h,k),\Phi^3_J\rangle
$$
where

$$B_J(h,k) :=
\sum_{I: |I|<|J|}
\langle h, \Phi^2_I\rangle
\langle k, \Phi^3_I\rangle
\Phi^1_I.
$$
The above formula is of type

$$\sum_{J\in \J}
\frac{1}{|J|^{1/2}}
a^1_J a^2_J a^3_J
$$
where $\J$ is an arbitray collection of dyadic intervals and such expressions can be estimated
by un upper bound of the form (see \cite{c})

\begin{equation}\label{se}
\prod_{i=1}^3
[\size_{\J}((a^i_J)_J)]^{1-\theta_i}\cdot
[\energy_{\J}((a^i_J)_J)]^{\theta_i}
\end{equation}
for any $0\leq \theta_1,\theta_2, \theta_3 < 1$ so that $\theta_1+\theta_2+\theta_3 = 1$
with the implicit constants depending on these ``theta parameters''.

The definitions of these ``sizes'' and ``energies'' are as follows:

$$\size_{\J}((a^i_J)_J) := \sup_{J\in\J} \frac{|a^i_J|}{|J|^{1/2}}$$
if $(a^i_J)_J$ is of ``$\Phi$ type'' (meaning that the corresponding implicit $\Phi^i_J$ functions are of ``$\Phi$ type''), and

$$\size_{\J}((a^i_J)_J) := \sup_{J\in\J}
\frac{1}{|J|}
\left\|\left(
\sum_{J'\in\J; J'\subseteq J}
\frac{|a^i_{J'}|^2}{|J'|}\chi_{J'}
\right)^{1/2}
\right\|_{1,\infty}
$$
if $(a^i_J)_J$ is of ``$\Psi$ type''.

Also, the ``energy'' is defined by

$$\energy_{\J}((a^i_J)_J) := \sup_{n\in\Z} \sup_{\D}2^n (\sum_{J\in\D}|J|)$$
where $\D$ either ranges over those collections of disjoint dyadic intervals $J\in\J$
for which

$$\frac{|a^i_J|}{|J|^{1/2}} \geq 2^n$$
in the ``$\Phi$ case'', or it ranges over the collection of disjoint dyadic intervals $J\in\J$
having the property that

$$
\frac{1}{|J|}
\left\|\left(
\sum_{J'\in\J; J'\subseteq J}
\frac{|a^i_{J'}|^2}{|J'|}\chi_{J'}
\right)^{1/2}
\right\|_{1,\infty}\geq 2^n
$$
in the ``$\Psi$ case''.

In the case of $(\langle f,\Phi^1_J\rangle)_J$ or $(\langle g,\Phi^2_J\rangle)_J$
sequences, there are ways to estimate further these sizes and energies, either by certain
averages of $f$ and $g$ or by the $L^1$-norms of $f$ and $g$ \cite{c}. The case of
$(\langle B_J(h,k),\Phi^3_J\rangle)_J$ is more complicated, since the inner function
depends on the interval $J$. The main observation here (called ``the bi-est trick'' in \cite{mtt:walshbiest}, \cite{mtt:fourierbiest})
is that this dependence can actually be factored out. 

More precisely, assume that one wants to estimate the size of such a sequence, and that the
suppremum is attained for an interval $J_0$. The size then becomes

\begin{equation}\label{26}
\frac{1}{|J_0|}
\left\|\left(
\sum_{J\in\J; J\subseteq J_0}
\frac{| \langle B_J(h,k),\Phi^3_J\rangle)  |^2}{|J|}\chi_{J}
\right)^{1/2}
\right\|_{1,\infty}.
\end{equation}
From the definition of $B_J(h,k)$ we see that terms of type 
$\langle \Phi^1_I, \Phi^3_J\rangle$ for $|I|<|J|$ are implicit in the above expression.
By using Plancherel, one then sees that one must have $\omega^3_J\subseteq \omega^1_I$ for such
a term to be nonzero. Denote now by $\I_0$ the set of all dyadic intervals $I\in\I$
for which there exists $J\subseteq J_0$ so that $\omega^3_J\subseteq \omega^1_I$
\footnote{These ``omega intervals''are the ``frequency intervals'' which support the Fourier 
transform of the corresponding functions. }.
We then observe that

$$\langle B_J(h,k),\Phi^3_J\rangle = \langle B_{\I_0}(h,k),\Phi^3_J\rangle$$
for any $J\subseteq J_0$ where

$$
B_{\I_0}(h,k) :=
\sum_{I\in\I_0}
\frac{1}{|I|^{1/2}}
\langle h,\Phi^2_I\rangle
\langle k,\Phi^3_I\rangle
\Phi^1_I.
$$
The reason for this is that if a pair of type $\langle \Phi^1_I, \Phi^3_J\rangle$
appears for which one has the opposite inclusion $\omega^1_I\subseteq \omega^3_J$, 
for some $I\in\I_0$ then, by the definition of $\I_0$ there exists another interval
$J'\subseteq J_0$ so that $\omega^3_{J'}\subseteq \omega^1_I\subseteq \omega^3_J$.
But then, $\omega^3_{J'}$ and $\omega^3_J$ would be one strictly inside the other,
which contradicts the frequency structure of the $(\omega^3_J)_J$ intervals.

This means that (\ref{26}) equals to

$$
\frac{1}{|J_0|}
\left\|\left(
\sum_{J\in\J; J\subseteq J_0}
\frac{| \langle B_{\I_0}(h,k),\Phi^3_J\rangle)  |^2}{|J|}\chi_{J}
\right)^{1/2}
\right\|_{1,\infty}
$$
which can be estimated in terms of a certain average of $B_{\I_0}(h,k)$, which itself
after a duality argument can be further estimated by using again this time a ``local variant'' of
the general upper bound (\ref{se}). And then, a similar reasoning (based on the ``bi-est trick'') helps to understand the energies.

After that, to estimate the other 4-linear form $\Lambda^{\#}$ corresponding to (\ref{25}), one applies again
the same generic estimate (\ref{se}) but this time the ``biest trick'' is no longer effective
and some other ``ad hoc'' arguments are necessary. The point here is that all these forms
can indeed be estimated with upper bounds which are independent on $\#$, which makes the
whole sum over $\#$ convergent in the end. For more details, see the original paper \cite{c}.

We would like to end the article with the observation that the flag paraproduct which naturally appears in the study of the 2D quadratic NLS
satisfies the same $L^p$ estimates as the operator $T_{ab}$ in Theorem \ref{main}. More precisely, we have

\begin{theorem}\label{nls}
The flag paraproduct $T(f_1, f_2, f_3)$ defined by

$$T(f_1, f_2, f_3)(x):=
\int_{\R^3}a(\xi_1, \xi_2) b(\xi_1, \xi_2, \xi_3)
\widehat{f_1}(\xi_1)
\widehat{f_2}(\xi_2)
\widehat{f_3}(\xi_3)
e^{2\pi i x(\xi_1 +\xi_2 +\xi_3)}
d \xi
$$
maps $L^{p_1}\times L^{p_2}\times L^{p_3} \rightarrow L^p$ boundedly, for every $1<p_1, p_2, p_3 <\infty$ with $1/p_1+1/p_2+p_3 = 1/p$.
\end{theorem}

\begin{proof}
Before starting the actual proof, we should mention that the $\R^d$ variant of the theorem holds also true and that this extension to euclidean spaces of arbitrary dimension
is really straightforward. We shall describe the argument in the one dimensional case, to be consistent with the rest of the paper. The main point is to simply realize that the 
discrete model operators studied in \cite{c} are enough to cover this case also. Let us assume that the symbols $a(\xi_1, \xi_2)$ and $b(\xi_1,\xi_2,\xi_3)$
are given by

$$a(\xi_1, \xi_2) = \sum_{k_1}
\widehat{\Phi_{k_1}}(\xi_1)
\widehat{\Phi_{k_1}}(\xi_2)$$
and

$$b(\xi_1,\xi_2,\xi_3) = \sum_{k_2}
\widehat{\Phi_{k_2}}(\xi_1)
\widehat{\Phi_{k_2}}(\xi_2)
\widehat{\Phi_{k_2}}(\xi_3)$$
as usual \footnote{Again, as we pointed out earlier, modulo some minor technical issues, one can always assume that this is the case.}.

As a consequence, we have

\begin{equation}\label{27}
a(\xi_1, \xi_2) b(\xi_1,\xi_2,\xi_3) = \sum_{k_1, k_2}
\widehat{\Phi_{k_1}}(\xi_1)
\widehat{\Phi_{k_1}}(\xi_2)
\widehat{\Phi_{k_2}}(\xi_1)
\widehat{\Phi_{k_2}}(\xi_2)
\widehat{\Phi_{k_2}}(\xi_3).
\end{equation}
The are two cases:

\underline{Case 1: $k_2<<k_1$.}

This case is simple since we must either have $k_1 \sim k_2$ (and the corresponding expression generates a classical paraproduct)
or both $(\widehat{\Phi_{k_1}}(\xi_1))_{k_1}$ and $(\widehat{\Phi_{k_1}}(\xi_2))_{k_1}$ are of 
``$\Phi$ type'' which is impossible.

\underline{Case 2: $k_1<<k_2$.}

In this case we have that both $(\widehat{\Phi_{k_2}}(\xi_1))_{k_2}$ and $(\widehat{\Phi_{k_2}}(\xi_2))_{k_2}$ are of ``$\Phi$ type'' which implies
in particular that $(\widehat{\Phi_{k_2}}(\xi_3))_{k_2}$ must be of ``$\Psi$ type''. We then rewrite (\ref{27}) as 

\begin{equation}\label{28}
a(\xi_1, \xi_2) b(\xi_1,\xi_2,\xi_3) = \sum_{k_1<< k_2}
\widehat{\Phi_{k_1}}(\xi_1)
\widehat{\Phi_{k_1}}(\xi_2)
\widehat{\Phi_{k_2}}(\xi_1)
\widehat{\Phi_{k_2}}(\xi_2)
\widehat{\Phi_{k_2}}(\xi_3).
\end{equation}
As in the ``$a(\xi_1,\xi_2)b(\xi_2,\xi_3)$ case'' we would have liked instead of (\ref{28}) to face an expression of type

\begin{equation}\label{29}
\sum_{k_1<< k_2}
\widehat{\Phi_{k_1}}(\xi_1)
\widehat{\Phi_{k_1}}(\xi_2)
\widehat{\Phi_{k_2}}(\xi_1+\xi_2)
\widehat{\widetilde{\Phi_{k_2}}}(\xi_1+\xi_2)
\widehat{\Phi_{k_2}}(\xi_3).
\end{equation}
Indeed, as before, if we faced this instead, we would have rewritten it as

\begin{equation}\label{30}
\sum_{k_1<< k_2}
\widehat{\Phi_{k_1}}(\xi_1)
\widehat{\Phi_{k_1}}(\xi_2)
\widehat{\widetilde{\widetilde{\Phi_{k_1}}}}(\xi_1+\xi_2)
\widehat{\widetilde{\widetilde{\widetilde{\Phi_{k_2}}}}}(\xi_1+\xi_2)
\widehat{\Phi_{k_2}}(\xi_3)
\end{equation}
and this, as we have seen, generates discrete models of type (\ref{24}) which have been understood in \cite{c}.

We show now how can one transform (\ref{28}) into an expression closer to (\ref{29}). The idea is once again based on using a careful Taylor decomposition,
this time for the functions $\widehat{\Phi_{k_2}}(\xi_1)$ and $\widehat{\Phi_{k_2}}(\xi_2)$.
We write

$$\widehat{\Phi_{k_2}}(\xi_1) = \widehat{\Phi_{k_2}}(\xi_1+ \xi_2) +
\frac{\widehat{\Phi_{k_2}}'(\xi_1+ \xi_2)}{1!}(-\xi_2) + 
\frac{\widehat{\Phi_{k_2}}''(\xi_1+ \xi_2)}{2!}(-\xi_2)^2 + ... +
\frac{\widehat{\Phi_{k_2}}^M(\xi_1+ \xi_2)}{M!}(-\xi_2)^M +
R_M(\xi_1, \xi_2)$$
and similarly

$$\widehat{\Phi_{k_2}}(\xi_2) = \widehat{\Phi_{k_2}}(\xi_1+ \xi_2) +
\frac{\widehat{\Phi_{k_2}}'(\xi_1+ \xi_2)}{1!}(-\xi_1) + 
\frac{\widehat{\Phi_{k_2}}''(\xi_1+ \xi_2)}{2!}(-\xi_1)^2 + ... +
\frac{\widehat{\Phi_{k_2}}^M(\xi_1+ \xi_2)}{M!}(-\xi_1)^M +
\widetilde{R_M}(\xi_1, \xi_2).$$
Now, if we insert these two formulae into (\ref{28}), we obtain (most of the time) expressions whose general terms are of type

$$
\widehat{\Phi_{k_1}}(\xi_1)
\widehat{\Phi_{k_1}}(\xi_2)
\frac{\widehat{\Phi_{k_2}}^l(\xi_1+ \xi_2)}{l!}(-\xi_2)^l
\frac{\widehat{\Phi_{k_2}}^{\widetilde{l}}(\xi_1+ \xi_2)}{\widetilde{l}!}(-\xi_1)^{\widetilde{l}}
\widehat{\Phi_{k_2}}(\xi_3) :=
$$

$$
\frac{1}{2^{k_2 l}}
\frac{1}{2^{k_2\widetilde{l}}}
\widehat{\Phi_{k_1}}(\xi_1)(-\xi_1)^{\widetilde{l}}
\widehat{\Phi_{k_1}}(\xi_2)(-\xi_2)^l
\widehat{\Phi_{k_2, l}}(\xi_1+\xi_2)
\widehat{\Phi_{k_2, \widetilde{l}}}(\xi_1+\xi_2)
\widehat{\Phi_{k_2}}(\xi_3) :=
$$

$$
\frac{2^{k_1 l}}{2^{k_2 l}}
\frac{2^{k_1\widetilde{l}}}{2^{k_2\widetilde{l}}}
\widehat{\Phi_{k_1, \widetilde{l}}}(\xi_1)
\widehat{\Phi_{k_1, l}}(\xi_2)
\widehat{\Phi_{k_2, l, \widetilde{l}}}(\xi_1+\xi_2)
\widehat{\Phi_{k_2}}(\xi_3)
$$
for $0\leq l,\widetilde{l}\leq M$.
Since $k_1<<k_2$, we can assume as before that $k_2 = k_1 + \#$, with $\#$ greater or equal than (say) $1000$.
In particular, if at least one of $l, \widetilde{l}$ is strictly bigger than zero, the corresponding operator can be written as

$$\sum_{\#=1000}^{\infty}
2^{-\#(l+\widetilde{l})} T_{\#}^{l,\widetilde{l}}.$$
Now, as we pointed out before, the analysis of each $T_{\#}^{l,\widetilde{l}}$ can be reduced to the analysis of the model operators (\ref{25}) which have been understood in 
\cite{c}. They are all bounded, with upper bounds which are uniform in $\#$ and this allows one to simply sum the above implicit geometric sum.
In the case when both $l$ and $\widetilde{l}$ are equal to zero, then the corresponding symbol is precisely of the form (\ref{29}) and the operator
generated by it can be reduced (as we already mentioned) to the model operators (\ref{24}) of \cite{c}.

Finally, the operators whose symbols are obtained when at least one of $R_M(\xi_1,\xi_2)$ or $\widetilde{R_M}(\xi_1,\xi_2)$ enters the picture, can be 
all estimated by the Coifman - Meyer theorem.

\end{proof}


\begin{thebibliography}{99}

\bibitem{ablowitzsegur}Ablowitz M., Segur H., {\it Solitons and the inverse scattering transform},
SIAM Studies in Mathematics, [1981].





\bibitem{meyerc}Coifman R. R., Meyer, Y.,
{\it Op\'erateurs multilin\'eaire}, Hermann, Paris, [1991].





\bibitem{ck1}Christ M., Kiselev A., {\it WKB asymptotic behaviour of almost all generalized eigenfunctions of one dimensional Schr\"{o}dinger operators},
J. Funct. Anal., vol. 179, 426-447, [2001].



\bibitem{ck2}Christ M., Kiselev A., {\it Maximal functions associated to filtrations},
J. Funct. Anal., vol. 179, 409-425, [2001].


\bibitem{gms}Germain P., Masmoudi N., Shatah J., {\it Global solutions for 3D quadratic Schr\"{o}dinger equations}, preprint [2008].




\bibitem{gt}Grafakos L., Torres R.,
{\it Multi-linear Calder\'{o}n-Zygmund theory}, Adv. Math., vol. 165,
124-164, [2002].

\bibitem{kp}Kato T., Ponce G., {\it Commutator estimates and the Euler and Navier-Stokes equations},
Comm. Pure Appl. Math., vol. 41, 891-907, [1988].


\bibitem{ks}Kenig C., Stein E.,
{\it Multilinear estimates and fractional integration}, Math. Res. Lett.,
vol. 6, 1-15, [1999].




\bibitem{laceyt1}Lacey M., Thiele C., {\it 
$L^p$ estimates on the bilinear Hilbert transform for $2<p<\infty$},
Ann. of Math., vol. 146, pp. 693-724, [1997].


\bibitem{laceyt2}Lacey M., Thiele C., {\it 
On Calderon's conjecture},
Ann. of Math., vol. 149, pp. 475-496, [1999].

\bibitem{c}Muscalu C., {\it Paraproducts with flag singularities I. A case study},
Rev. Mat. Iberoamericana, vol. 23 705-742, [2007].

\bibitem{c1}Muscalu C., {\it Paraproducts with flag singularities II. The general case},
in preparation.

\bibitem{cct}Muscalu C., Tao T., Thiele C., {\it Multilinear operators
given by singular symbols}, J. Amer. Math. Soc., vol. 15, 469-496, [2002].

\bibitem{mtt:walshbiest}Muscalu C., Tao T., Thiele C., 
{\it $L^p$ estimates for the biest I.  The Walsh case}, Math. Ann., vol. 329,
401-426, [2004].

\bibitem{mtt:fourierbiest}Muscalu C., Tao T., Thiele C., 
{\it $L^p$ estimates for the biest II.  The Fourier case}, Math. Ann., 
vol. 329, 427-461, [2004].



\bibitem{mptt:biparameter}Muscalu C., Pipher J., Tao T., Thiele C.,
{\it Bi-parameter paraproducts}, Acta Math., vol. 193, 269-296, [2004].

\bibitem{mptt:multiparameter}Muscalu C., Pipher J., Tao T., Thiele C.,
{\it Multi-parameter paraproducts}, Revista Mat. Iberoamericana, vol 22, 963-976, [2006].



\bibitem{mtt:multiest}Muscalu C., Tao T., Thiele C.,
{\it Multi-linear multipliers associated to simplexes of arbitrary length}, preprint [2007].



\bibitem{nrs}Nagel A., Ricci F., Stein E., {\it Singular integrals with flag kernels on quadratic CR manifolds}, J. Funct. Anal.. vol. 181, 129-181, [2001].

\bibitem{ns}Nagel A., Stein E., {\it The $\overline{\partial}_b$ - complex on decoupled boundaries in $\C^n$}, Ann. of Math. (2), vol. 164, 649-713, [2006].

\bibitem{simon}Simon B., {\it Bounded eigenfunctions and absolutely continuous spectra for one-dimensional Schr\"{o}dinger operators},
Proc. Amer. Math. Soc., vol. 124, 3361-3369, [1996].


\bibitem{stein}Stein E.,
{\it Harmonic analysis: real-variable methods, orthogonality and
oscillatory integrals}, Princeton University Press, [1993].



\end{thebibliography}
\end{document}